\documentclass[letterpaper,12pt]{article}
\usepackage[hideboxes]{boxedeps}

\usepackage{latexsym}
\newcommand{\pic}[2]{\BoxedEPSF{#1 scaled #2}}

\def\CA{{\cal A}}
\def\CC{{\cal C}}
\def\CS{{\cal S}}

\def\Xor{\pic{Xor.ART} {500}}
\def\Yor{\pic{Yor.ART} {500}}
\def\Ior{\pic{Ior.ART} {500}}
\def\Sigmaior{\pic{Sigmaior.ART} {500}}
\def\Pij{\pic{Pij.ART} {500}}
\def\Wij{\pic{Wij.ART} {500}}

\def\Idor{\pic{Idor.ART} {500}}
\def\Rcurlor{\pic{Rcurlor.ART} {500}}

\def\Murphy{\pic{Murphy.ART} {500}}
\def\Murphysum{\pic{Murphysum.ART} {500}}
\def\Murphycurl{\pic{Murphycurl.ART} {500}}
\def\Murphysumless{\pic{Murphysumless.ART} {500}}

\def\Gnplus{\pic{Gnplus.ART} {500}}
\def\Gshift{\pic{Gshift.ART} {500}}
\def\Gcross{\pic{Gcross.ART} {500}}

\def\Aid{\pic{Aid.ART} {500}}
\def\Aa{\pic{Aa.ART} {500}}
\def\Aaa{\pic{Aaa.ART} {500}}
\def\Aainv{\pic{Aainv.ART} {500}}
\def\Annulus{\pic{Annulus.ART} {500}}
\def\Aeh{\pic{Aeh.ART} {500}}
\def\Ahe{\pic{Ahe.ART} {500}}
\def\Wn{\pic{Wn.ART} {500}}
\def\Vn{\pic{Vn.ART} {500}}
\def\WC{\pic{WC.ART} {500}}

\def\XC{\pic{XC.ART} {500}}
\def\XCYC{\pic{XCYC.ART} {500}}
\def\XYC{\pic{XYC.ART} {500}}
\def\TC{\pic{TC.ART} {500}}
\def\PhiX{\pic{PhiX.ART} {500}}
\def\PsiD{\pic{PsiD.ART} {500}}
\def\PsiX{\pic{PsiX.ART} {500}}
\def\PsiXY{\pic{PsiXY.ART} {500}}
\def\PsiXPsiY{\pic{PsiXPsiY.ART} {500}}

\newtheorem{theorem}{Theorem}[section]
\newtheorem{corollary}[theorem]{Corollary}
\newtheorem{lemma}[theorem]{Lemma}

\newenvironment{remark}{\par\smallskip%
\noindent\textbf{Remark.}\  }%
{\par\smallskip}

\newenvironment{Def}{\par\smallskip%
\noindent\textbf{Definition.}\  }%
{\par\smallskip}

\newenvironment{proof}[1][{}]{\par\smallskip%
\noindent\textit{Proof #1: }\  }
{\hfill$\Box$\par\smallskip}

\begin{document}
\begin{center}
      {\Large\bf Skein theory and the Murphy operators}

\bigskip
HUGH R. MORTON 

\medskip
 {\it Department of Mathematical Sciences, University of Liverpool,\\
     Peach St, Liverpool, L69 7ZL, England. }\\
{\tt morton@liv.ac.uk}
\end{center}
\vskip 1cm

\begin{abstract}
The Murphy operators in the Hecke algebra
$H_n$ of type $A$ are explicit commuting
elements whose sum generates the centre.
 They can be represented by simple tangles in the Homfly skein theory
version of $H_n$.
 In this paper I present a single tangle which represents their sum,
and which is obviously central. As a consequence
 it is possible to identify a natural basis for the Homfly skein of the
annulus, $\CC$.

Symmetric functions of the Murphy operators are also central in $H_n$.
I define geometrically a  homomorphism
from
$\CC$ to the centre of each algebra $H_n$, and
find an element in $\CC$, independent of $n$, 
whose image  is the $m$th power sum of the
Murphy operators. Generating function
techniques are used to describe images of
other elements of
$\CC$ in terms of the Murphy operators, and to
demonstrate  relations among other natural
skein elements.

\phantom{some stuff to separate the keywords from the abstract}

\end{abstract}

\section*{Introduction.}

The Hecke algebra $H_n$ of type $A$ is  a deformed version of
the group algebra ${\bf C}[S_n]$ of the symmetric
group. It has a simple
skein theory model in terms of $n$-tangles and
the Homfly skein relations.

Certain sums of transpositions,
\[m(j)=\sum_{i=1}^{j-1} (i\,j) \in
{\bf C}[S_n],\ j=2,\ldots,n,\]
 known as Jucys-Murphy elements,
appeared in the work of Jucys \cite{Jucys} in
1974, and later in work of Murphy \cite{Murphy}. 
These elements all commute, and every symmetric
polynomial in them can be shown to lie in the
{\em centre} of the algebra ${\bf C}[S_n]$. For
example
$m(2)=(1\,2), m(3)=(1\,3)+(2\,3)$ and
\[m(2)m(3)=(1\,3\,2)+(1\,2\,3)=m(3)m(2).\]

Dipper and James \cite{Dipper} used a simple deformation of the
transpositions  to define analogous
elements $M(j)\in H_n$, which they called the
{\em Murphy operators}. These elements again all
commute, and symmetric polynomials in them belong
to the centre of $H_n$. Dipper and James
 showed that for generic values of
the deformation parameter these account for the
whole of the centre. 
Katriel, Abdessalam and Chakrabarti \cite{Katriel}
noted that in fact any central element can be
expressed as a polynomial in just the sum
$M=\sum_{j=2}^n M(j)$ of the Murphy operators.
This is a stronger result than for the
non-deformed algebra ${\bf C}[S_n]$, although
the centre has the same dimension
in each case, given by $\pi(n)$, the number of
partitions of $n$.

In this paper  I present  a skein theory
version of the Murphy operators and  their sum, finishing with an
elegant representation of their power sums and other symmetric
polynomials.  Starting in section 2 with a choice of
$n$-string braids corresponding to the
transpositions I  exhibit a  braid $T(j)$
representing each of the individual Murphy
operators $M(j)$, following Ram \cite{Ram}, and then
  a very natural simple
$n$-tangle $T^{(n)}$ which  represents their sum $M$, up to
a linear combination with the identity element in
$H_n$ in each case. From the tangle viewpoint it becomes immediately
clear that $M$ is central.

The Homfly skein of the
annulus, $\CC$, has been exploited for many
years, for example in the work of Jun Murakami,
and myself and Short, in constructing and analysing
Homfly-based  invariants  by the use
of satellites. These invariants cover the same ground
as the quantum $sl(N)$ invariants for all choices
of module and have been developed variously in
algebraic and skein theoretic ways, for example
by Wenzl \cite{Wenzl}, Aiston and myself \cite{AistonMorton} and
 Blanchet
\cite{Blanchet}. 

 In section 3 I  show how the tangle view of the
Murphy operators gives a direct means for identifying the best
linear basis for $\CC$ when
studying these quantum invariants and establishing
its multiplicative properties. I then introduce a rather overlooked
homomorphism from $\CC$ to the centre of $H_n$. Using generating
function methods I show how to realise other symmetric functions of
the Murphy operators as images of explicit elements of $\CC$. The
principal result is the identification of an element $P_m\in\CC$,
independent of $n$, whose image determines the $m$th power sum of the
Murphy operators in $H_n$.

Section 4 gives the final details of the results of section 3, based
on a further underused skein, $\CA$, whose unexpected algebraic
properties allow for some satisfyingly clean proofs.

\section{The skein models.}

The skein theory model of $H_n$ is
based on the framed Homfly skein relations, in
their simplest form
\begin{eqnarray*}
\Xor\  -\ \Yor& =&(s-s^{-1})\ \Ior\\
\mbox{and }\  \Rcurlor &=&v^{-1}\ \Idor.
\end{eqnarray*}

The Homfly skein ${\cal S}(F)$ of a planar
surface $F$, with some designated input and
output boundary points, is defined as linear
combinations of oriented tangles in $F$, modulo
these two local relations, and Reidemeister moves
II and III. The coefficient ring can be taken
as $\Lambda={\bf Z}[v^{\pm1}, s^{\pm1}]$ with
powers of $s^k-s^{-k}$ in the denominators.

Every skein admits a {\em mirror map},
$\overline{\phantom{w}}:\CS(F)\to\CS(F)$ induced by switching all
crossings in a tangle, coupled with inverting $v$ and $s$ in $\Lambda$.

Write $R_n^n(v,s)$ for the skein ${\cal S}(F)$ of
$n$-tangles,  where $F$
is a rectangle with $n$ inputs at the bottom and
$n$ outputs at the top. Composing $n$-tangles
induces a product which makes $R_n^n(v,s)$ into
an algebra. It has a linear basis of $n!$
elements, and is isomorphic to the Hecke algebra
$H_n(z)$. This algebra has a presentation with
generators
$\{\sigma_i\}, i=1,\ldots,n-1$ satisfying the
braid relations
\begin{eqnarray*}
\sigma_i\sigma_j&=&\sigma_j\sigma_i,
\quad |i-j|>1,\\
\sigma_i\sigma_{i+1}\sigma_i&=&\sigma_{i+1}\sigma_i\sigma_{i+1},
\end{eqnarray*} and the quadratic relations
$\sigma_i^2=z\sigma_i+1$.

In $R_n^n$ the generators are the elementary
braids
\[\sigma_i\ =\ \Sigmaior\]
and the parameter $z$ is $s-s^{-1}$, giving the
alternative form $(\sigma_i-s)(\sigma_i+s^{-1})=0$
for the quadratic relations.

In the special case $z=0$ the Hecke algebra
reduces to ${\bf C}[S_n]$, with $\sigma_i$
becoming the transposition
$(i\;i+1)$. For general $i,j$
we can view the transposition $(i\,j)$
 diagrammatically \[(i\,j)\ =\ \Pij\ .\] The corresponding
positive permutation braid \[\omega_{(i\,j)}\ =\ \Wij\] with
all crossings positive realises this
transposition, as in \cite{AistonMorton}. Regarded via
$R_n^n$  as an element of $H_n(z)$, it
becomes the transposition $(i\,j)$ when we set
$z=0$.

\begin{Def} The {\em Murphy operator} $M(j)\in R^n_n,\ j=2,\ldots,n$ is
defined as
\[M(j)=\sum_{i=1}^{j-1}\omega_{(i\,j)}.\]
\end{Def}
These elements certainly project to the
Jucys-Murphy elements $m(j)$ in ${\bf C}[S_n]$.
 With a
bit of algebraic work they  can be shown to
commute, and their sum
\[M=\sum_{j=2}^n M(j)\quad
=\sum_{i<j}\omega_{(i\,j)}
\] can be shown to lie in the centre of $R^n_n$. These facts are
immediately obvious from the skein representatives $T(j)$ and
$T^{(n)}$, discussed shortly. 

\subsection{Variants of the Hecke algebras.}

A simple adjustment of the skein relations, as in \cite{AistonMorton},
allows for a skein model of $H_n$ whose parameters can be readily
adapted to match any of the different appearances of the algebra.

Extend the coefficient ring to include an invertible parameter $x$,
and define the skein $R^n_n(x,v,s)$ by linear combinations of oriented
$n$-tangles modulo the relations
\begin{eqnarray*}
x^{-1}\Xor\  -\ x\Yor& =&(s-s^{-1})\ \Ior\\
\mbox{and }\  \Rcurlor &=&xv^{-1}\ \Idor.
\end{eqnarray*}

There is a natural algebra homomorphism $R_n^n(v,s) \to R^n_n(x,v,s)$
 induced by replacing each tangle $T$
in the skein 
$R_n^n(v,s)$ by $x^{wr(T)}T$ in $R^n_n(x,v,s)$, where $wr(T)$ is the
writhe of $T$.

We then get corresponding versions of the Murphy elements in our new
skein; for example the new version of $T(j)$ is represented by
$x^{2j-2}T(j)$ in $R^n_n(x,v,s)$.

The parameters in $R^n_n(x,v,s)$ can be adjusted to match the exact
version of $H_n$ under study. Taking $x=s,q=s^2, v=1$ we
get the common algebraist's version $H_n(q)$ of \cite{Dipper}, while
the endomorphism rings arising from the fundamental representation of
$sl(N)_q$ and its standard $R$-matrix correspond to the choice
$x=s^{-1/N}, v=s^{-N}$ with $s=e^{h/2}$.

The choice $x=v$ gives the Homfly skein relation which is invariant
under Reidemeister move I, and $x=q^r, s=q, v=a$ gives the version
used by Kawagoe in \cite{Kawagoe}.

 The elements
$M(j)$ above do indeed
  become the Murphy operators used by 
Dipper and James under the isomorphism of
$R^n_n(v,s)$ with the algebraic version $H_n(q)$ of the Hecke algebra
where $x=s,q=s^2, v=1$.

\section{Geometric views of the Murphy operators.}
The Murphy operators $M(j)$ and their sum $M$ have simple skein
representatives as shown here,  using $R^n_n(v,s)$ as the model for
the Hecke algebra $H_n$.

\begin{theorem} \label{Murphy} {\rm (Ram)}. The Murphy operator
$M(j)$ can be represented in $H_n$ by a single
braid $T(j)$, up to linear combination with the
identity.
\end{theorem}

\begin{theorem} \label{Sum}  {\rm (Morton)}. The sum $M$ of the Murphy operators can be
represented in $H_n$ by a single tangle $T^{(n)}$, again up to linear
combination with the identity.

\end{theorem}

\begin{proof} [ of theorem \ref{Murphy}] In \cite{Ram} Ram notes  that $M(j)$ can
{\em almost} be represented by the single braid
$T(j)$ pictured, where `almost' means that $M(j)$
is a linear combination of $T(j)$ and the
identity.

Set \[T(j)\ =\ \Murphy \] in $H_n$, for each $j\le n$ including the
case
$T(1)=1$. Skein theory shows quickly that $T(j)-1=zM(j)$, giving
\[M(j)=\frac{T(j)-1}{s-s^{-1}}.\]
\end{proof}

The elements $T(j)$ will do equally well in place
of $M(j)$, so long as $z\neq0$, in other words
away from ${\bf C}[S_n]$.

 It is clear
geometrically that the elements $T(j)$ all
commute.  It is not immediately clear
that the {\em sum} of these elements is in the centre
of $H_n$, although their product is the full
twist, a well-known central element.

\begin{proof}  [of theorem \ref{Sum}]
Set \[T^{(n)}=\Murphysum \]in $R_n^n\cong H_n(z)$, with the coefficient
ring extended to include $v^{\pm 1}$.

Apply the skein relation at one crossing to get
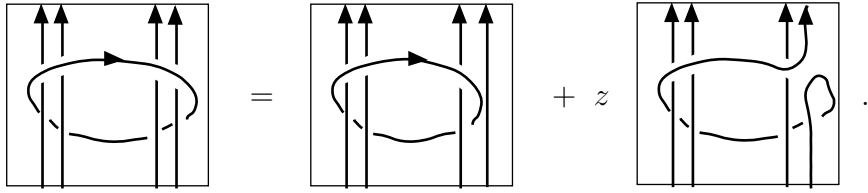
\[ \Murphysum\ =\ \Murphysumless \ + \ z\ \Murphycurl.\] Hence
$T^{(n)}=T^{(n-1)}+zv^{-1}T(n)$, using the
standard inclusion of $H_{n-1}$ in $H_n$ to interpret $T^{(n-1)}$ as
an element of $H_n$. Then
\[T^{(n)}-T^{(0)} =zv^{-1}\sum_{j=1}^nT(j),\]
 by induction on $n$. Here $T^{(0)}$ is just a scalar multiple of the
identity, represented by a single disjoint simple loop alongside the
identity braid. Since a trivial loop in any Homfly skein contributes the
scalar $\delta=(v^{-1}-v)/z$ we can rewrite the equation as 
\[T^{(n)}-\frac{v^{-1}-v}{z}=zv^{-1}\sum_{j=1}^n
T(j)=z^2v^{-1}M+nzv^{-1}.\]

Then
\[T^{(n)}=z^2v^{-1}M+\left(nzv^{-1}+\frac{v^{-1}-v}{z}\right)\times
1,\] and hence, again away from $z=0$, we can write $M$ as a linear
combination of
$T^{(n)}$ and $1$.
\end{proof}

From this representation it is 
quite obvious that $T^{(n)}$, and thus $M$,  is in the
centre of
$H_n$.

There is a known set of idempotent elements,
$E_\lambda$, one for each partition $\lambda$ of
$n$, which were originally described
algebraically by Gyoja \cite{Gyoja}. Skein pictures of
these based on the Young diagram for $\lambda$
were given by Aiston and myself \cite{AistonMorton}. We
 showed there, in a skein based counterpart of the results of
\cite{Katriel},  that
$ME_\lambda=m_\lambda E_\lambda$ for an explicit scalar $m_\lambda$
and that
$m_\lambda\ne m_\mu$ when $\lambda\ne\mu$. It
follows at once that
$T^{(n)}E_\lambda=t_\lambda E_\lambda$, where the
scalars $t_\lambda$ are different for each
partition $\lambda$. 

\section{The skein of the annulus.}

The Homfly
skein of the annulus, ${\cal C}$, as discussed in
\cite{Turkey} and originally in \cite{Turaev},  is defined as linear combinations of
diagrams in the annulus, modulo the Homfly skein relations. The
element $X\in\CC$ will be indicated on a diagram as
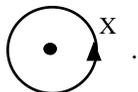
\[\XC\ .\]

 We shall
make use of  the results of section 2  to study $\CC$, and also to find
simple skein representatives for the power sums of the Murphy
operators in $H_n$ in a way which is independent of
$n$.

The best known relation of $\CC$ with the Hecke algebra $H_n$ is the
closure map
$R_n^n\to {\cal C}$, induced by taking a tangle $T$ to its
closure $\hat {T}$ in the annulus, defined by \[\hat {T}\ =\ \TC\ .\]
This is a linear map, whose image we will call ${\cal C}_n$.

The skein ${\cal C}$ has a product induced by
placing one annulus outside another,  under which $\cal C$ becomes a
commutative algebra;
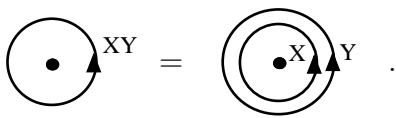
\[\XYC\ =\ \XCYC\ .\]
Write $A_m \in {\cal C}$ for the closure of
the
$m$-braid \[\sigma_{m-1}\cdots \sigma_2\sigma_1\ =\ \Gcross\ .\]
An explicit  spanning set of
${\cal C}_n$ consists of the monomials in 
$\{A_j\}$ of total weight $n$, where
$A_j$ has weight $j$. Such monomials
correspond bijectively with partitions of $n$.

A very simple skein theory construction
determines a natural linear map $\varphi:{\cal
C}\to{\cal C}$, induced by taking a diagram $X$
in the annulus and linking it once with a simple
loop to get
\[\varphi(X)\ =\ \PhiX\ .\]

It is instructive to look at the eigenvectors of
$\varphi$.  Clearly $\varphi$ carries ${\cal
C}_n$ into itself. What is more, if we take an
element $S\in H_n$ with closure $\hat{S}\in {\cal
C}_n$ and compose it with the central element $T^{(n)}$ then
$ST^{(n)}$ has closure $\varphi(\hat{S})$. 

\begin{theorem} \label{Phi} The eigenvalues of $\varphi|{\cal C}_n$ are
all distinct.
\end{theorem}

\begin{proof}Take
$S=E_\lambda$ to see that the closure of $T^{(n)}E_\lambda$ is
$t_\lambda\hat{E_\lambda}$ and also
$\varphi(\hat{E_\lambda})$. The element
$Q_\lambda=\hat{E_\lambda}\in{\cal C}_n$ is then
an eigenvector of $\varphi$ with eigenvalue
$t_\lambda$. Now there are $\pi(n)$ of these eigenvectors,  (the number
of partitions of $n$), and
the eigenvalues are all distinct, by \cite{AistonMorton}. Since 
${\cal C}_n$ is spanned by $\pi(n)$ elements we
can deduce that the elements $Q_\lambda$ form a
basis for ${\cal C}_n$, and that the eigenspaces of $\varphi$ are all
$1$-dimensional.
\end{proof}

 It follows that {\em any} element of ${\cal C}_n$
which is an eigenvector of $\varphi$ must be a multiple of
some $Q_\lambda$.

The elements $h_i\in {\cal C}$, discussed below, arise
 as
$Q_\lambda$ for the Young diagram $\lambda$ consisting of a
single row with $i$ cells.  Kawagoe \cite{Kawagoe} constructed an
element
$s_\lambda$ of
${\cal C}_n$ for each partition $\lambda$ of $n$,
as a Schur polynomial of the sequence $\{h_i\}$, in Macdonald's context
of symmetric functions, \cite{Macdonald}. He then used skein theory
to show that each
$s_\lambda$ is also an eigenvector of the map
$\varphi$. 
His elements
$s_\lambda$ can now be identified almost at once with the elements
$Q_\lambda$, by theorem \ref{Phi}.

This gives an immediate proof of the multiplicative properties of the
elements $Q_\lambda$ in the skein of the annulus, since the elements
$s_\lambda$ automatically multiply according to the Littlewood-Richardson
rules for Young diagrams \cite{Macdonald}. Conversely it guarantees
that Kawagoe's elements
$s_\lambda$ are eigenvectors for a  wide range of skein maps on ${\cal
C}$ as well as $\varphi$, using the skein properties of $E_\lambda$
from
\cite{AistonMorton}. 

  Recently Lukac, \cite{Lukac}, 
 has  been able to simplify the proof
that $s_\lambda$ is an eigenvector of $\varphi$, by using algebraic
properties of the skein $\CA$ of the annulus with an input on one
boundary and an output on the other. 
Discussion and use of the skein $\CA$ in section \ref{A} 
leads to simple descriptions for other symmetric functions of the
Murphy operators, derived from a direct interplay between the
skein of the annulus and the centre of the Hecke algebras.

\subsection{Symmetric functions and the skein of the annulus.}

In this section I recall some explicit results about elements in the
Hecke algebras and their closure in $\CC$, and their interpretation in
the context of symmetric functions, following the methods of
Macdonald. This leads to a simple description in terms of generating
functions of some interrelations in $\CC$, and also of the elements to
be used later to represent the power sums of the Murphy operators.

The starting point is the description of the two simplest idempotents
in $H_n$, corresponding to the single row and column Young diagrams.
These are given algebraically in \cite{Jones}; here I  use the
skein version as described in \cite{Turkey} in terms of the positive
permutation braids $\omega_\pi,\ \pi\in S_n$. Define two
quasi-idempotents by
\[a_n=\sum_{\pi\in S_n}s^{l(\pi)}\omega_\pi, \quad b_n =\sum_{\pi\in
S_n}(-s)^{-l(\pi)}\omega_\pi,\] where $l(\pi)=wr(\omega_{\pi})$ is the
writhe of the braid $\omega_{\pi}$. 

The following result is straightforward. 
\begin{lemma} We can write \[a_n=a_{n-1}\gamma_n,\] where
$\gamma_n=1+s\sigma_{n-1}+s^2\sigma_{n-1}\sigma_{n-2}
+\cdots+s^{n-1}\sigma_{n-1}\cdots\sigma_1$. 
\end{lemma} 
We have $\gamma_{n+1}=1+s\sigma_n\gamma_n$, and also an immediate skein
relation
\[ \Gnplus\ =\ \Gshift\ +\ s^{n}\ \Gcross\ ,\] in $R^{n+1}_{n+1}$.

The next lemma is proved in \cite{AistonMorton}.
\begin{lemma} \label{Rowidem} For any braid
$\beta\in B_n$ we have $a_n\beta=\varphi_s(\beta)a_n=\beta a_n$,  where
$\varphi_s(\beta)=s^{wr(\beta)}$.
\end{lemma}
Similar results, with $s$ replaced by $-s^{-1}$, hold for $b_n$.

The element $a_n$ then satisfies
\[a_n^2=\varphi_s(a_n)a_n=\varphi_s(a_{n-1})\varphi_s(\gamma_n)a_{n}.\]

Since $\varphi_s(\gamma_n)=1+s^2+\cdots+s^{2n-2}=s^{n-1}[n]$, we have
an immediate corollary.
\begin{corollary} \label{Idempotent} We can write
\[s^{n-1}[n]h_n=h_{n-1}\gamma_n,\] where $h_n=a_n/\varphi_s(a_n)$ is
the true idempotent.
\end{corollary}

The element $h_n\in H_n$ is the idempotent which corresponds to the
single row with $n$ cells. The single column idempotent is given from
$h_n$ by using $-s^{-1}$ in place of $s$.

With a slight abuse of notation  write $h_n\in {\cal C}$
in place of $\hat{h}_n$ for the closure of this element in ${\cal C}$,
and $e_n\in{\cal C}$ for the closure of the single column
idempotent.

\begin{remark} Aiston uses  the notation $Q_{c_n}$ and
$Q_{d_n}$  in \cite{Aiston}, in place of $e_n$ and $h_n$ which are
used here to suggest
 the terminology and techniques of symmetric functions from
\cite{Macdonald}.
\end{remark}

Write \[H(t)=1+\sum_{n=1}^\infty h_n t^n\] for the generating function
of the elements $\{h_n\}$, regarded as a formal power series with
coefficients in ${\cal C}$, and similarly set
$E(t)=1+\sum_{n=1}^\infty e_n t^n$. In \cite{Aiston} Aiston showed that
\begin{theorem} \[E(-t)H(t)=1,\] as power series in ${\cal C}$.
\end{theorem}

If we regard the elements $h_n$ formally as the $n$th complete
symmetric functions in a large number $N$ of variables $x_1,
\ldots,x_N$, setting
\[H(t)=\prod_{i=1}^N \frac{1}{1-x_it}\ ,\] then $E(t)=\prod(1+x_it)$
and
$e_n$ is the $n$th elementary symmetric function in the variables.

Following Macdonald, the Schur function $s_\lambda(x_1,\ldots,x_n)$ for each
partition $\lambda$
can be expressed as a polynomial in the elements $h_n$ 
 which does not depend on $N$ for large enough
$N$, and is therefore determined formally by the series $H(t)$ and
$\lambda$ as an element $s_\lambda$ of ${\cal C}$. This is the element
$s_\lambda$ used by Kawagoe, which is identified  with the
idempotent closures $Q_\lambda$ by Lukac in \cite{Lukac}.

The $\pi(n)$ elements $\{s_\lambda\}$ with $|\lambda|=n$ form a linear
basis for ${\cal C}_n$, and can be expressed in terms of the monomials
of weight $n$ in $\{h_r\}$ by the classical Jacobi-Trudy formulae
\cite{Macdonald}. 

Indeed the skein ${\cal C}^+$, defined as $\cup{\cal C}_n$, is spanned
by all monomials in $\{h_r\}$, and can be interpreted as the ring of
symmetric polynomials in an unlimited number of variables $\{x_i\}$.

Monomials in the geometrically simpler closed braid elements
$\{A_m\}$ also span
${\cal C}_n$, and it is  thus interesting to relate these directly to
$\{h_r\}$. An attractive formula connecting these two generating sets
can be derived from corollary \ref{Idempotent}. 

\begin{theorem} \label{AH}
Write $A(t)=1+z\sum_{m=1}^\infty A_mt^m$, with $z=s-s^{-1}$.

Then \[A(t)=\frac{H(st)}{H(s^{-1}t)}=H(st)E(-s^{-1}t).\]
\end{theorem} 
The proof will be given in section \ref{A}.

\begin{lemma} The elements $\{h_n\}$ are invariant under
the mirror map switching crossings and inverting $s$ and $v$.
\end{lemma}
\begin{proof}In $H_n$ we have
$\sigma_ia_n=\varphi_s(\sigma_i)a_n=sa_n$. Then
$\overline{\sigma_ia_n}=\sigma_i^{-1}\overline{a}_n=s^{-1}\overline{a}_n$
giving $\sigma_i\overline{a}_n=s\overline{a}_n$. It follows that
$\beta\overline{a}_n=\varphi_s(\beta)\overline{a}_n$ for any $\beta$
and so $a_n\overline{a}_n=\varphi_s(a_n)\overline{a}_n$.

On the other hand $a_n\overline{a}_n=a_n\varphi_s(\overline{a}_n)$ by
lemma \ref{Rowidem}, so that
\[h_n=\frac{a_n}{\varphi_s(a_n)}
=\frac{\overline{a}_n}{\varphi_s(\overline{a}_n)}=\overline{h}_n.\]
\end{proof}
\begin{corollary} The inverse series to $A(t)=1+z\sum A_mt^m$ is
$\overline{A}(t)=1-z\sum \overline{A}_mt^m$, where $\overline {A}_m$
is the  negative closed braid corresponding to $A_m$.
\end{corollary}
\begin{proof} Apply the mirror map in theorem \ref{AH} to get
\[\overline{A}(t)={\overline{H(st)}}\ /\ {\overline{H(s^{-1}t)}}
={H(s^{-1}t)}\ /\ {H(st)}=A(t)^{-1}.\]
\end{proof}

In our later description of the power sums of Murphy operators we use
the element
$P_m
\in {\cal C}$ which can be interpreted formally as the $m$th power sum,
$P_m=\sum x_i^m$, of the variables. This is determined unambiguously as
a polynomial in
$\{h_r\}$, and hence an element of ${\cal C}$, by Newton's power series
equation
\[\sum_{m=1}^\infty \frac{P_m}{m}t^m=\ln H(t).\]
Aiston showed in \cite{Adams} that $[m]P_m$ is the sum of $m$ closed
$m$-string braids given  by switching the first
$i$ crossings of the closed braid $A_m$ from positive to negative, for
$i=0,\ldots,m-1$.  Her proof depends on some quantum group translations
to identify
$s_\lambda$ and $Q_\lambda$, and to write $A_m$ in terms of
$s_\lambda$. In another article I give a direct skein theory proof of
her result, using theorem \ref{AH}.

\subsection{Symmetric functions of the Murphy operators.}
The  main result in this paper, besides theorem \ref{Sum} and its
consequences, is an expression for the
$m$th power sum of the Murphy operators in $H_n$ in terms of the
element
$P_m\in {\cal C}$.

For this I use a previously unremarked relation between the
Hecke algebras and the skein of the annulus.
This relation takes the form of a very natural homomorphism $\psi_n$
from
${\cal C}
$ to the centre of each algebra $H_n$.  

 The diagram
\[D\ =\ \PsiD\] determines a  map $\psi_n:{\cal C}\to H_n$,
induced by placing $X\in \CC$  around the circle in
$D$ and the identity of $H_n$ on the arc, to get
\[\psi_n(X)\ =\ \PsiX\ \in H_n.\]
 It is clear that 
\[\psi_n(XY)\ =\ \PsiXY\ =\ \PsiXPsiY\ =\psi_n(X)\psi_n(Y),\]
so that
$\psi_n$ is an algebra homomorphism.
Moreover the elements $\psi_n(X)$ obviously all lie in the
centre of $H_n$, and in fact it follows from theorem \ref{Power}
below that the image of $\psi_n$ consists af all symmetric polynomials
in the Murphy operators and so, by \cite{Dipper}, makes up the whole of
the centre in the generic case.

We already know from theorem \ref{Sum} that the sum $\sum T(j)$ is
essentially $T^{(n)}$, and we can write $T^{(n)}=\psi(X)$ with $X=h_1$,
represented by a single string around the core of the annulus. We may
then look for elements
$X_2$ with
$\psi(X_2)=\sum T(j)^2$ and more generally $\psi(X_m)=\sum
T(j)^m$.

These can certainly be found for each fixed
$n$, using for  example a suitable polynomial in $h_1$ which depends on
$n$. Much more surprising is the result of theorem \ref{Power},
that there exists an element
$X_m$ {\em independent of $n$} with the property
that
\[\psi_n(X_m)=\sum_{j=1}^n T(j)^m,\] up to a
multiple of the identity. Consequently any other
symmetric function of the Murphy operators can be
achieved as $\psi(X)$ where the choice of $X$ is
essentially independent of $n$.

In fact, we can take $X_m=P_m$, to get the following explicit result,
proved in section \ref{A}.
\begin{theorem} \label{Power} For any $n$ we have
\[\psi_n(P_m) -\psi_0(P_m) = (s^m-s^{-m})v^{-m}\sum_{j=1}^n T(j)^m.\]
\end{theorem}

This is particularly satisfactory in that $P_m$ itself arises as a
power sum, although not of any identifiable objects in the ring $\cal
C$ itself.

 Since any element
$X\in {\cal C}^+$ can be written as a polynomial in
$\{P_m\}$ we can then express $\psi(X)$ as a symmetric function of the
Murphy operators $T(j)$.  It is useful to be able to write the
elements 
 $\psi(h_i)$ in terms of the Murphy operators; theorem \ref{Power}
leads to a  compact expression using  formal power series with
coefficients in the centre of the Hecke algebra. 

I shall do this in terms of
the {\em Murphy series} $HM(t)$ of the Hecke algebra,
whose coefficients, lying in the centre of $H_n$,  are defined to be
the complete symmetric functions of the elements
$T(j)$. Thus the Murphy series in $H_n$ can be
written explicitly as \[HM(t)=\prod_{j=1}^n
(1-T(j)t)^{-1}.\] The series has a formal inverse, $EM(-t)$, where the
coefficients of
$EM(t)=\prod(1+T(j)t)$ are the elementary symmetric functions of the
elements $T(j)$, also in the centre of
$H_n$.

\begin{theorem}

\[\psi_n(H(t))=\psi_0(H(t))\frac{HM(sv^{-1}t)}{HM(s^{-1}v^{-1}t)}.\]

\end{theorem}

\begin{proof} It is enough to establish that the logarithms of the
two sides are equal.  Now  \[\ln
(\psi_n(H(t)))=\psi_n(\ln(H(t)))=\sum_{m=1}^\infty
\frac{\psi_n(P_m)}{m} t^m,\] while \[\ln (HM(t))=\sum_{m=1}^\infty
\frac{\sum T(j)^m}{m} t^m.\] The coefficients of $t^m$ in $\ln
(\psi_n(H(t)))$ and in
\begin{eqnarray*}
&&\ln\left(\psi_0(H(t))\frac{HM(sv^{-1}t)}{HM(s^{-1}v^{-1}t)}\right)\\
&=&\ln\left(\psi_0(H(t))\right)+
\ln\left(HM(sv^{-1}t)\right)-\ln\left(HM(s^{-1}v^{-1}t)\right)
\end{eqnarray*}
are equal, by theorem
\ref{Power}.
\end{proof}

\section{The annulus with two boundary points.} \label{A}

To prove theorems \ref{AH} and \ref{Power} I introduce a third skein 
$\CS(F)$, whose underlying surface is the annulus with a single input
on one boundary component and corresponding output on the other, as
shown,

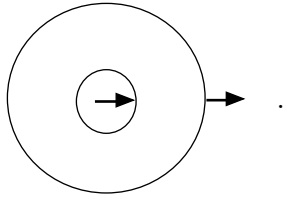
\begin{center}
\Annulus\ .
\end{center}

Write $\CA$ for this skein, whose elements are represented as linear
combinations of oriented tangles in the annulus consisting of a single
arc and a number of closed curves.  As for the skein $\CC$ of the
annulus there is a product on $\CA$ induced by placing one annulus
outside the other.  The identity element $e\in \CA$ is represented by
the tangle \[ \Aid\ .\]
Write $a\in \CA$ for the element represented by \[ \Aa\ .\] From this
we can construct $a^m$ for each $m\in {\bf Z}$, giving for example
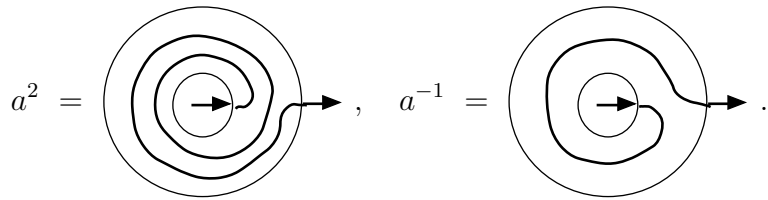
\[ a^2\ =\ \Aaa\ ,\quad a^{-1}\ =\ \Aainv\ .\]

There are {\em two} bilinear products, $l:\CC\times\CA\to\CA$ and
$r:\CA\times\CC\to\CA$ induced by placing an element of $\CC$
respectively under or over an element of $\CA$. For example, $h_1\in
\CC$, represented by a single counterclockwise loop, gives
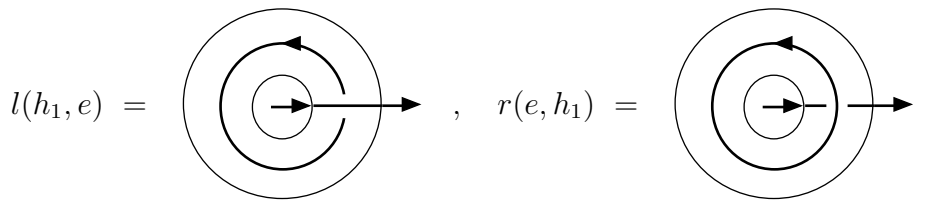
\[l(h_1,e)\ =\ \Ahe\ ,\quad r(e,h_1)\ =\ \Aeh\ .\]

From the skein relation in $\CA$ we have
$l(h_1,e)-r(e,h_1)=(s-s^{-1})a$.

\subsection{Algebraic properties of $\CA$.}
Kawagoe \cite{Kawagoe} and other authors have used a skein which is
linearly isomorphic to $\CA$, based on
the annulus with input and output on the {\em same} component. The
significant advantage of $\CA$ lies in its algebraic properties.

Both $l$ and  $r$ are algebra homomorphisms, since it is clear
diagrammatically that $l(c_1,a_1)l(c_2,a_2)=l(c_1c_2,a_1a_2)$, while
$r$ behaves similarly.

\begin{theorem} The algebra $\CA$ is commutative.
\end{theorem}

\begin{proof}
Unlike the case of $\CC$ this is not immediately clear. Using standard
skein theory techniques we can represent any element of $\CA$ as a
linear combination of tangles consisting of a totally descending arc
lying over a number of closed curves. Each such tangle represents
$l(c_m,a^m)=l(c_m,e)a^m$ for some $m$ and some $c_m\in \CC$. The
general element of $\CA$ can then be written as a Laurent polynomial
\[\sum_{m\in {\bf Z}}l(c_m,e)a^m\]
in $a$, with coefficients in the commutative subalgebra
$l(\CC,e)\subset \CA$. Since $a$ commutes with $l(\CC,e)$ it follows
that any two elements of $\CA$ commute.
\end{proof}

The subalgebras $l(\CC,e)$ and $r(e,\CC)$ are both isomorphic to
$\CC$, but they are not equal. Their difference determines a sort of
commutator map $[\ ,e]:\CC\to\CA$, defined by $[c,e]=l(c,e)-r(e,c)$.

We have already noted that $[h_1,e]=(s-s^{-1})a$. 
We can deduce
theorem \ref{Power}
about the power sums of Murphy operators from the next theorem.
\begin{theorem}\label{Powercomm} For $m\ge 1$ we have
$[P_m,e]=(s^m-s^{-m})a^m$.
\end{theorem}

\subsection{Skein interaction between $H_n$, $\CA$ and $\CC$.}
We can make use of wiring diagrams, as described in
\cite{AistonMorton}, to induce linear maps between skeins. A {\em
wiring diagram} $W$ is an inclusion of one surface $F$ into another
$F'$, along with a fixed diagram of curves in $F'$ which avoid $F$, and
connect any distinguished input and output points on the boundaries of
both, respecting orientation.  A wiring $W$ induces a linear map,
denoted $\CS(W):\CS(F)\to \CS(F')$, using $W$ to extend any diagram in
$F$ to a diagram in $F'$.

To avoid over-elaborate notation  I shall write $W$ in place of
$\CS(W)$ where there is no risk of confusion.

The diagram \[W_n\ =\ \Wn\] with $n$ strings running around the
annulus then induces a linear map $W_n:R^{n+1}_{n+1}\to \CA$, by
inserting an $(n+1)$-tangle in the empty box.

When an element $S\in H_n$ is included in $H_{n+1}$ in the standard
way, using their skein versions, we can see from a diagram that
$W_n(S)=W_{n-1}(S)a$.
Then $W_n(1)=a^n$, where $a\in \CA$ is the element shown above.

\begin{theorem}
The elements $W_n(h_{n+1}), W_n(h_n)$ and $l(h_n,e)$ in $\CA$ satisfy
the linear relation
\[[n+1]W_n(h_{n+1})=s^{-1}[n]W_n(h_n)+l(h_n,e).\]
\end{theorem}
\begin{proof}
The relation \[ \Gnplus\ =\ \Gshift\ +\ s^{n}\ \Gcross\] gives
$W_n(h_n\gamma_{n+1})=W_n(\gamma_n h_n)+s^nW_n(h_n
\sigma_n\cdots\sigma_1)$,  after taking $\gamma_n$ from the second
diagram around the strings of the wiring.

Now $\gamma_n h_n=s^{n-1}[n]h_n$ and
$W_n(h_n\sigma_n\cdots\sigma_1)=l(h_n,e)$. Combined with the equation
$s^n[n+1]h_{n+1}=h_n\gamma_{n+1}$ from corollary \ref{Idempotent} the
result follows at once.
\end{proof}

Now write $Y_n=[n+1]W_n(h_{n+1})$ and set $Y(t)=\sum_{n=0}^\infty
Y_nt^n$.
\begin{corollary}\label{HY} As power series with coefficients in $\CA$
we have
\[l(H(t),e)=(e-s^{-1}at)Y(t).\]
\end{corollary}
\begin{proof} Since $W_n(h_n)=W_{n-1}(h_n)a=aW_{n-1}(h_n)$ the relation
can be written $Y_n=s^{-1}aY_{n-1}+l(h_n,e)$. This gives
$Y(t)=s^{-1}atY(t)+l(H(t),e)$, and hence the result.
\end{proof}

 The mirror map, switching crossings and inverting $s$
and $v$, when applied to
$l(X,e)$ for any $X\in \CC$ gives $r(e,\overline X)$. Now
$\overline{H(t)}=H(t)$ and so 
$\overline{Y(t)}=Y(t)$, giving
\[r(e,H(t))=(e-sat)Y(t).\]

Subtracting this from the equation in corollary \ref{HY} then gives
$[H(t),e]=(s-s^{-1})atY(t)$, a result used by Lukac in \cite{Lukac}.

\begin{proof}[of theorem \ref{Powercomm}]
Take logarithms of the equations and then subtract. Then
\begin{eqnarray*}
\ln(e-s^{-1}at)-\ln(e-sat)&=&\ln\left(l(H(t),e)\right)
-\ln\left(r(e,H(t))\right)\\
&=& l(\ln(H(t)),e)-r(e,\ln(H(t)))\\
&=&\sum_{m=1}^\infty
[\frac{P_m}{m}t^m,e].
\end{eqnarray*}
Now $\ln(e-s^{-1}at)=-\sum_{m=1}^\infty s^{-m}a^mt^m/m$, and it only
remains to compare the coefficients of $t^m$ to see that
$[P_m,e]=(s^m-s^{-m})a^m$.
\end{proof}

\begin{proof}[of theorem \ref{Power}] Use the  wiring \[V_n\ =\ \Vn\]
with $n-1$ strings passing through the annulus to carry
$\CA$ into
$R^n_n$. Then for any $X$ in $\CC$ we have $V_n(l(X,e))=\psi_n(X)$ and
$V_n(r(e,X))=\psi_{n-1}(X)$ with the standard inclusion.
We can also see that $V_n(a)=v^{-1}T(n)$, and extend this to
$V_n(a^m)=v^{-m}T(n)^m$. Then
$\psi_n(P_m)-\psi_{n-1}(P_m)=(s^m-s^{-m})v^{-m}T(n)^m$.

Induction on $n$ completes the proof that
\[\psi_n(P_m)-\psi_0(P_m)=(s^m-s^{-m})v^{-m}\sum_{j=1}^n T(j)^m.\]
\end{proof}

To conclude this paper I give the proof of theorem \ref{AH}, using one
final wiring, this time to induce a map from $\CA$ to $\CC$.

Write
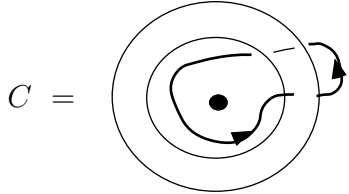
\[C\ =\ \WC\] for the wiring diagram shown, which provides a form of
closure for an element of $\CA$, and determines a linear map
$C:\CA\to\CC$. It is clear that $C(a^m) =A_{m+1}$ for every $m\ge 0$.
We can also see readily that  \[C(F\ r(e,X))=C(r(F,X))=C(F)X\] for
every $F\in \CA$ and $X\in \CC$.

\begin{proof} [of theorem \ref{AH}] The relation
$(e-sat)Y(t)=r\left(e,H(t)\right)$ gives
\[Y(t)=\left(\sum_{m=0}^\infty s^ma^mt^m\right)r\left(e,H(t)\right).\]
Then
\begin{eqnarray*} C\left(Y(t)\right)&=&C\left(\sum_{m=0}^\infty
s^ma^mt^m\right)H(t)\\
&=&\left(\sum_{m=0}^\infty s^mA_{m+1}t^m\right)H(t),
\end{eqnarray*}
giving
\[zst  C\left(Y(t)\right)=\left(A(st)-1\right)H(t).\]

Now $C\left(W_n(h_{n+1})\right)$ is the closure of the $(n+1)$-tangle
$\sigma_1\sigma_2\cdots\sigma_nh_{n+1}=s^nh_{n+1}$. Then
$C(Y_n)=[n+1]s^nh_{n+1}\in\CC$. It follows that
\begin{eqnarray*}
zst 
C\left(Y(t)\right)&=&\sum_{n=0}^\infty(s^{n+1}-s^{-(n+1)})s^{n+1}h_{n+1}t^{n+1}\\
&=&\sum_{n=0}^\infty(s^{2n+2}-1)h_{n+1}t^{n+1}\\
&=&H(s^2t)-H(t).
\end{eqnarray*}
This gives $H(s^2t)-H(t)=A(st)H(t)-H(t)$, and thus
$H(st)=A(t)H(s^{-1}t)$, to complete the proof.
\end{proof}
\medskip
\noindent
Version 1.2, February 2001.

\end{document}